\documentclass{article}

\usepackage{graphicx}
\usepackage{latexsym}
\usepackage[utf8]{inputenc}
\usepackage{amsfonts, amssymb}
\usepackage{amsmath}
\usepackage{enumerate}
\usepackage{graphicx}
\usepackage[labelsep=period,font=small]{caption}
\usepackage{mathrsfs}
\usepackage{authblk}
\usepackage{color}

\newtheorem{thm}{Theorem}[section]
\newtheorem{defn}[thm]{Definition}
\newtheorem{lemma}[thm]{Lemma}
\newtheorem{prop}[thm]{Proposition}

\newtheorem{rem}[thm]{Remark}
\newtheorem{cor}[thm]{Corollary}

\newcommand{\C}{\mathbb C}
\newcommand{\Om}{\Omega}
\newcommand{\om}{\omega}
\newcommand{\sm}{\setminus}
\newcommand{\ga}{\gamma}

\newcommand{\si}{\sigma}

\newcommand{\Laplace}{\Delta}
\newcommand{\BB}{\mathcal B}
\newcommand{\eps}{\epsilon}
\newcommand{\N}{\mathbb N}
\newcommand{\D}{\mathbb D}
\newcommand{\KK}{\mathcal K}
\newcommand{\de}{\delta}

\newcommand{\dH}{\ensuremath{{\operatorname{d_H}}}}
\renewcommand{\DH}{\ensuremath{{\operatorname{D_H}}}}
\renewcommand{\d}{\ensuremath{\operatorname{d}}}
\newcommand{\Reg}{\ensuremath{{\mathrm{Reg}}}}

\DeclareMathOperator{\Cpct}{Cap}
\DeclareMathOperator{\Co}{Co}
\begin{document}

\title{Julia Sets of Orthogonal Polynomials
\thanks{
The authors would like to thank
the Danish Council for Independent Research $|$ Natural Sciences for support via the
grant DFF -- 4181-00502.
The last author would also like to thank the Institute of Mathematical Sciences
of Stony Brook University for support and hosting during the 
writings of the paper.}
}

\author[Petersen, Pedersen, Henriksen and Christiansen]{Jacob Stordal Christiansen, Christian Henriksen, Henrik Laurberg Pedersen and Carsten Lunde Petersen
}

\date{\today}

\maketitle

\begin{abstract}

For a probability measure with compact and non-polar support in the complex plane we relate dynamical properties of the associated sequence of orthogonal polynomials
$\{P_n\}$ to properties of the support. More precisely we relate the Julia set of $P_n$ to the outer boundary of the support, the filled Julia set to the polynomial convex hull $K$ of the support, and the
Green's function associated with $P_n$ to the Green's function for the complement of $K$.

\end{abstract}



\section{Introduction and main results}
\label{intro}

In this paper, we study
orthonormal polynomials
$\{P_n(z)\}\equiv \{P_n(\mu; z)\}$ given by a Borel probability measure $\mu$ on $\C$ with 
compact and non-polar support $S(\mu)$. We relate the non-escaping set for $P_n$, the locus of non-normality (the boundary of the non-escaping set) for $P_n$,
and an associated Green's function to the support of the measure, 
getting a fairly complete picture of the limiting behavior of these objects as $n \to \infty$.

We build on the classical monograph \cite{StahlandTotik} by Stahl and Totik,
where the authors relate potential and measure theoretic properties of, e.g.,
the asymptotic zero distribution for the sequence of orthonormal polynomials
defined by $\mu$ to the potential and measure theoretic properties of the support of $\mu$.
We shall also use \cite{Randsford} as a reference to the basic concepts of potential theory in the complex plane.

Recall that $\{P_n(z)\}$ is the unique orthonormal sequence in $L^2(\mu)$ with
\begin{equation}
 \label{normalpol}
P_n(z)
= \ga_n z^n + \textrm{~lower order terms}, 
\end{equation}
where $\ga_n>0$.
\begin{defn}
Let $\BB$ denote the set of Borel probability measures on $\C$ with compact, non-polar
support. Furthermore, let $\BB_+\subset\BB$ be defined as
\begin{equation}
\BB_+ := \{\mu\in\BB \,|\, \limsup_{n\to\infty}\,\ga_n^{1/n} < \infty \},
\end{equation}
where $\ga_n$ is given in \eqref{normalpol}.
\end{defn}

For $\mu \in \BB$ we denote by
$\Om$ the unbounded connected component of $\C\sm S(\mu)$
and define
$$
K=\C\sm\Om,\quad   J=\partial K.
$$
The set $K$ is the \emph{filled} $S(\mu)$ and $J= \partial\Om \subset S(\mu)$ is the outer boundary of $S(\mu)$.
We shall also say that $S(\mu)$ is \emph{full} if $\C\sm S(\mu)$
has no bounded connected components.


Furthermore, we define the exceptional subset $E\subset S(\mu)$ by
\begin{equation}
E = \{ z \in J \,|\, z \textrm{ is \emph{not} a Dirichlet regular boundary point}\}.
\end{equation}
This set is an $F_\si$ polar subset, see \cite[Theorem 4.2.5]{Randsford}.
We let ${g_\Om}: \C \to [0,\infty)$ be the Green's function
for $\Om$ with pole at infinity (in short, just the Green's function for $\Om$).
This is the unique non-negative subharmonic function which is harmonic and positive
on $\Om$, zero precisely on $K\sm E$, (see \cite[Theorem 4.4.9]{Randsford}) 
and which satisfies
\begin{equation}
g_\Om(z) = \log|z| + O(1) \;\mbox{ at infinity}.
\end{equation}
Finally, we denote by $\om_J$ the equilibrium measure on $J$, which equals harmonic measure on
$\Om$ from $\infty$ and which is the distributional Laplacian $\Laplace g_\Om$
of the Green's function $g_\Om$.

We shall also use (see \cite[Section 1.2]{StahlandTotik}) the extended notion of the Green's function $g_B: \C\to [0,\infty)$
for an arbitrary connected Borel set $B\subset\C$ with bounded complement $L$ of positive logarithmic capacity, $\Cpct(L)>0$.
This is the unique non-negative subharmonic function which is harmonic and positive on the interior of $B$, satisfies
\begin{equation}
g_B(z) = \log|z| - \log\Cpct(L) + o(1)\; \textrm{ at infinity},
\end{equation}
and equals zero qu.\,e.~on $\C\sm B$. Here, qu.\,e.~is short for quasi everywhere meaning except on a polar set (\cite{Randsford} uses n.\,e., nearly every\-where).

Furthermore, for $\mu\in\BB$ we denote by $g_\mu: \C \to [0,\infty)$
{the minimal carrier Green's function for} $\mu$ (see \cite[Definition 1.1.1 and Lemma 1.2.4]{StahlandTotik}),
\begin{equation}
g_\mu(z) = \log|z| - \log c_\mu + o(1)\; \textrm{ at infinity},	
\end{equation}
where $c_\mu$ is the minimal carrier capacity.
Moreover, we denote by $E_\mu$ the exceptional set for $g_\mu$ defined by
\begin{equation}
E_\mu = \{ z\in S(\mu) \,|\, g_\mu(z) > 0 \}.
\end{equation}

The following fundamental result concerning the distribution of zeros of the orthogonal polymonials was originally obtained by Fej\'er in \cite{Fejer}; see also \cite[Lemma 1.1.3]{StahlandTotik}.
\begin{thm}
\label{thm:fejer}
If $\mu\in\BB$, then all zeros of the orthonormal polynomials $P_n$ 
are contained in the convex hull $\Co(S(\mu))$. Moreover, for any compact subset
$V\subset\Om$ the number of zeros of $P_n$ 
in $V$ is bounded as $n\to\infty$.
\end{thm}

Our main result, Theorem \ref{THMlimsupliminf}, concerns measures in the class $\BB_+$ and it is proved in Section \ref{sec:3}.
The first part of the theorem should be compared with \cite[Theorem 1.1.4]{StahlandTotik}, 
while the second part does not have an immediate counterpart in the classical theory of orthogonal polynomials.
We remark that $\BB_+$ is a large subclass of $\BB$ since only measures in $\BB$ with zero carrier capacity are left out.

Before stating our main result, some more notation is needed. We denote by $\Om_n$ the attracted basin of $\infty$ for $P_n$,
by $K_n = \C\sm \Om_n$ the filled Julia set, and by $J_n =\partial K_n = \partial\Om_n$
the Julia set.
Theorem \ref{THMlimsupliminf} is motivated by the following questions: What is the relation between $K$ and limits of $K_n$
and, similarly, what is the relation between $J$ and limits of $J_n$?
Inspired by \cite{Douady}, we answer these questions in terms of limits involving the Hausdorff distance on the space of compact subsets of $\C$
(see the beginning of Section \ref{sec:3} for details and the notions
of $\liminf$ and $\limsup$ of sequences of compact sets).
\begin{thm}
 \label{THMlimsupliminf}
Let $\mu\in\BB_+$. Then the following assertions hold.
\begin{enumerate}
 \item[(i)]
We have\begin{equation}
\limsup_{n\to\infty} K_n \subseteq \Co(K).
\end{equation}
Moreover, for any $\eps >0$ and with $V_\eps := \{z\in \C \,|\, g_\Om(z) \geq \eps \}$,
\begin{equation}
\lim_{n\to\infty}\Cpct(V_\eps\cap K_n) = 0.
\end{equation}
\item[(ii)] We have
\begin{equation}
\overline{J\sm E_\mu} \subseteq \liminf_{n\to\infty} J_n.
\end{equation}\end{enumerate}
\end{thm}




The figure below
illustrates Theorem \ref{THMlimsupliminf} 
in the case where $\mu$ is the equilibrium measure for the boundary of the boomerang-shaped white set $K$
in the top left image. The black fractal sets in the other images
are the Julia sets $J_{10}$, $J_{15}$, and $J_{20}$
(which in these cases appear to be equal to the filled Julia sets).
The Green's functions are visualised by
colouring alternating intervals of level sets blue and red.

We remark that equilibrium measures belong to a special class of measures,
the so-called {\it regular measures} to be discussed in Sections \ref{sec:4} and
\ref{sec:5}.

\begin{center}
\setlength\tabcolsep{1 pt}
\renewcommand{\arraystretch}{0.5}
\begin{tabular}{cc}
\includegraphics[width=56mm]{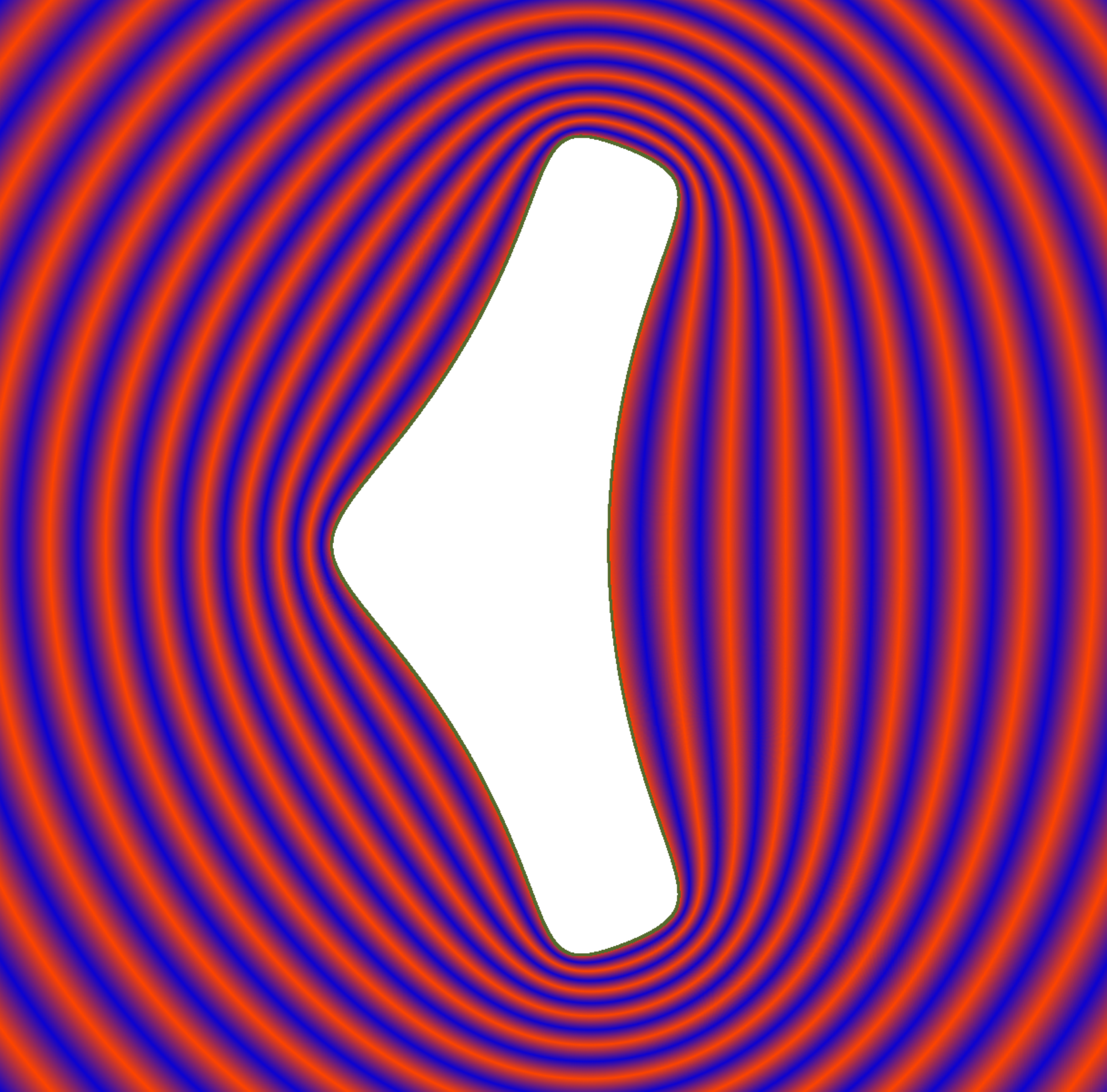}&
\includegraphics[width=56mm]{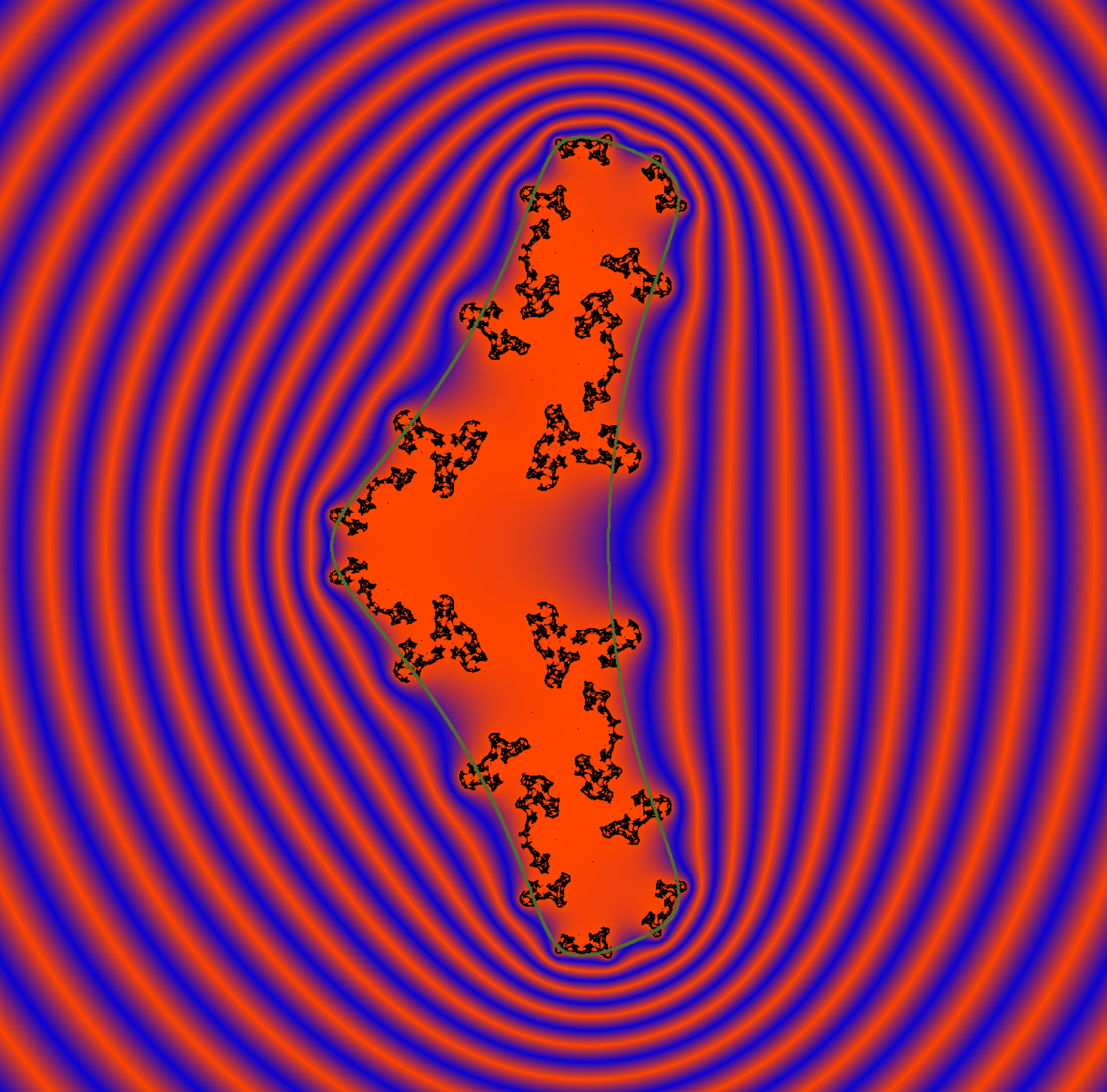}\\
\includegraphics[width=56mm]{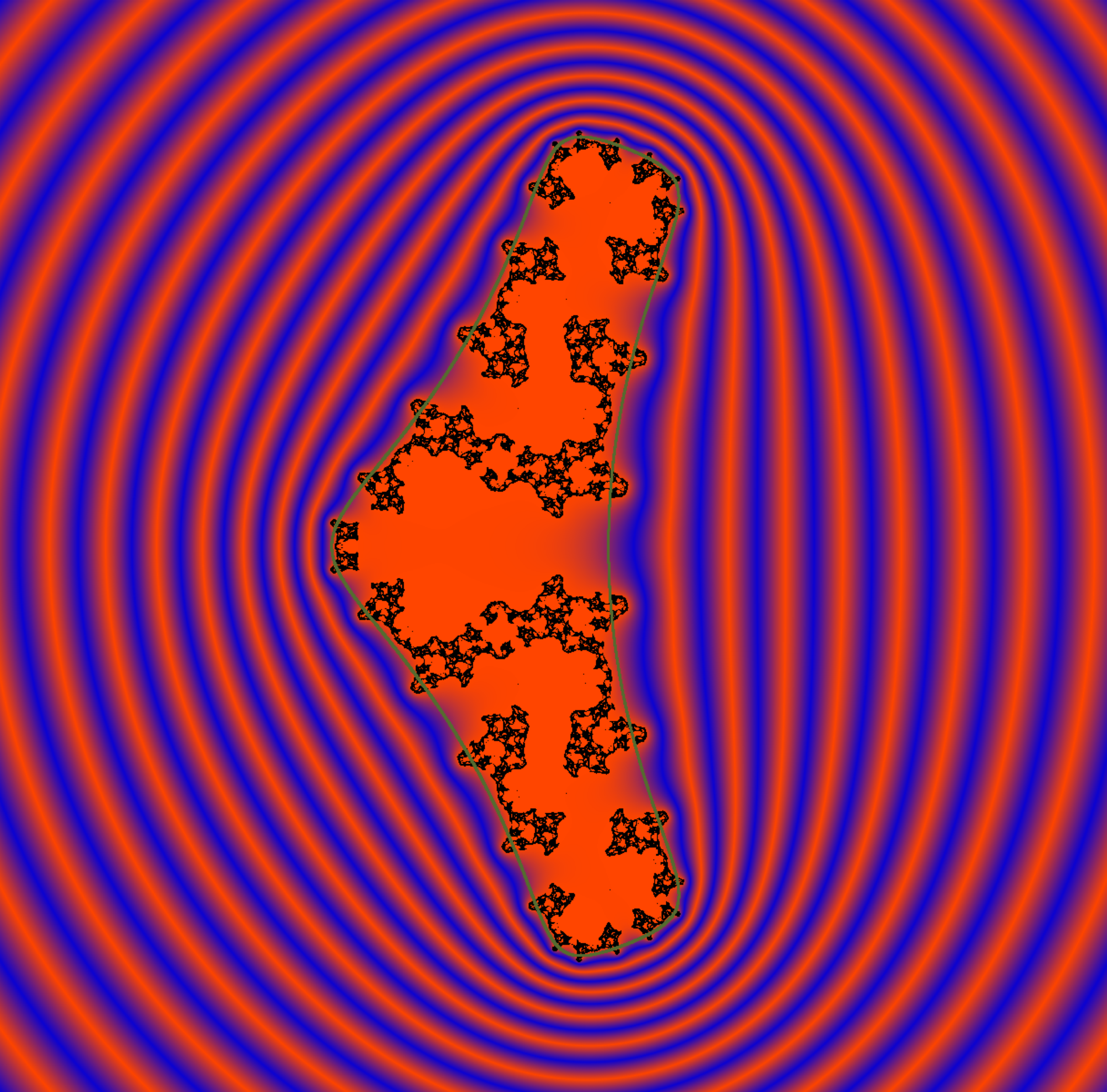}&
\includegraphics[width=56mm]{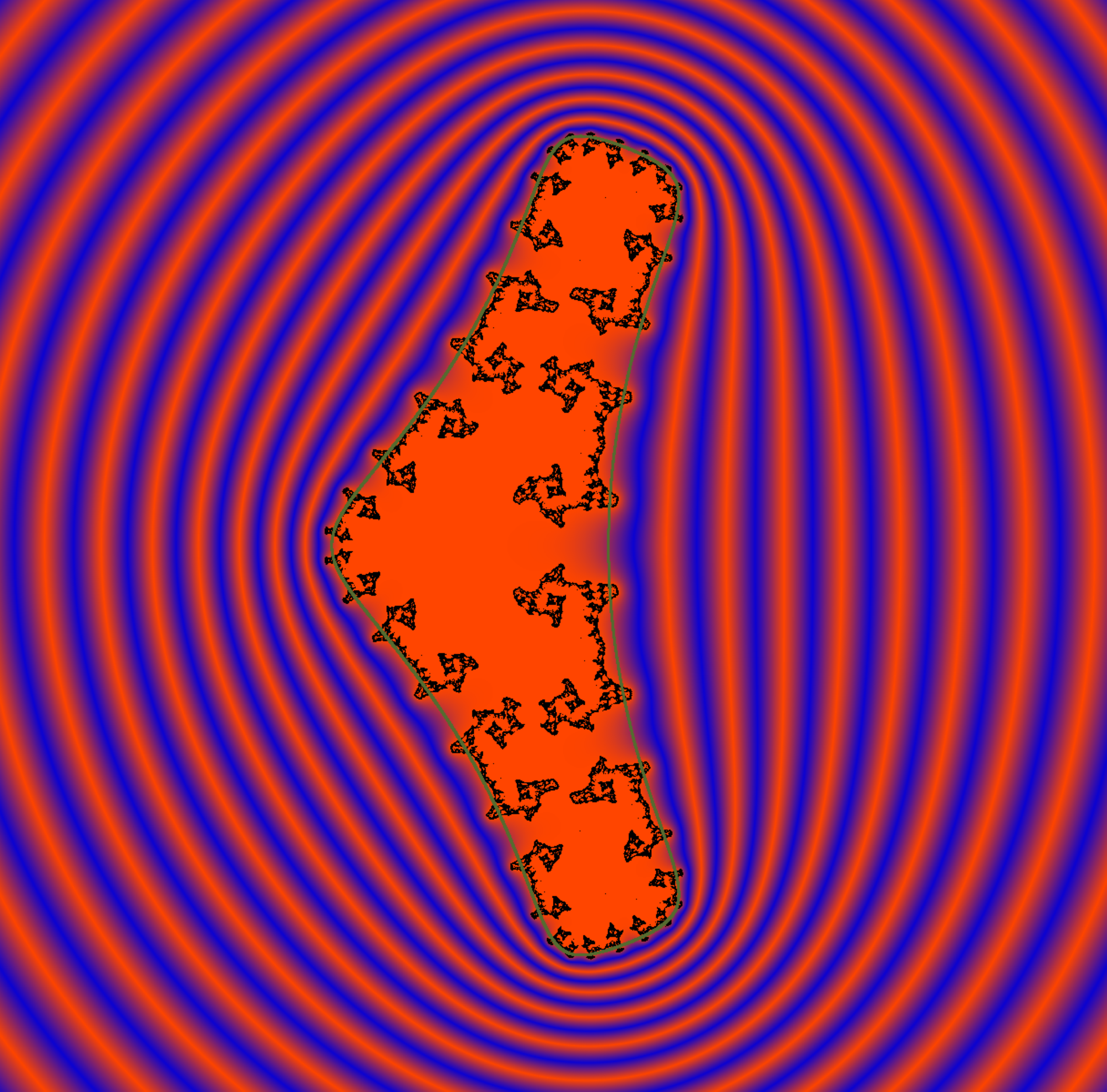}
\end{tabular}
\end{center}

\section{Polynomial dynamics and technical preparations} \label{sec:2}



For any polynomial $P$ of degree $d>1$, there clearly exists $R>0$ such that $|P(z)| \geq 2|z|$ for all $z$ with $|z|>R$.
Thus the orbit, $\{z_n\}$,  of such $z$ under iteration by $P$ converges to $\infty$.
The basin of attraction for $\infty$ for $P$, denoted $\Om_P$, may therefore be written as
\begin{equation}
\label{basinofinfty}
 \Om_P = \{z\in\C \,|\, P^k(z) \underset{k\to\infty}\longrightarrow \infty \}
= \bigcup_{k\geq 0} P^{-k}(\C\sm\overline{\D(0,R)}).
\end{equation}
Here $P^k = \overset{k \textrm{ times}}{\overbrace{P\circ P \circ \ldots \circ P}}$, whereas $P^{-k}$ denotes inverse image and $\D(0,R)$ is the open ball of radius $R$ centered at $0$.
It follows immediately that $\Om_P$ is open and completely invariant, that is,
$P^{-1}(\Om_P) = \Om_P = P(\Om_P)$.
Denote by $K_P = \C\sm\Om_P\subseteq\overline{\D(0,R)}$ the filled Julia set for $P$ and
by $J_P = \partial\Om_P = \partial K_P$ the Julia set for $P$.
Then $K_P$ and $J_P$ are compact and also completely invariant.
Clearly, any periodic point (i.e., a solution of the equation $P^k(z) = z$ for some $k\in\N$) belongs to $K_P$,
so that $K_P$ is non-empty.
It follows from \eqref{basinofinfty} that the filled Julia set $K_P$ can also be described as the nested intersection
\begin{equation}
\label{filledJuliaset}
K_P = \bigcap_{k\geq 0} P^{-k}(\overline{\D(0,R)}).
\end{equation}
To ease notation we denote the Green's function for $\Om_P$ with pole at infinity by $g_P$  (and not by $g_{\Om_P}$).
It follows from \eqref{filledJuliaset} that $g_P$ satisfies
\begin{equation}
g_P(z) = \lim_{k\to\infty} \frac{1}{d^k}\log^+(|P^k(z)|/R) = \lim_{k\to\infty} \frac{1}{d^k}\log^+|P^k(z)|.
\end{equation}
Here and elsewhere, $\log^+$ is the positive part of $\log$.
Thus $g_P$  vanishes precisely on $K_P$ 
and hence (\cite[Theorem 4.4.9]{Randsford}) every point in $J_P$ is a Dirichlet regular boundary point of $\Om_P$.
Moreover, denoting the leading coefficient of $P$ by $\gamma$,
\begin{equation}
g_P(P(z)) = d\cdot g_P(z)
\quad\textrm{and}\quad
\Cpct(K_P) = \frac{1}{|\ga|^{\frac{1}{d-1}}}.
\end{equation}
%
When $P=P_n$, we thus have in our notation
\begin{equation}\label{capacityformula}
\frac{1}{\ga_n^{\frac{1}{n-1}}} = \Cpct(K_n).
\end{equation}
As
\begin{equation} \label{infsup}
\liminf_{n\to\infty} \ga_n^{\frac{1}{n-1}} = \liminf_{n\to\infty} \ga_n^{\frac{1}{n}}
\quad\textrm{and}\quad
\limsup_{n\to\infty} \ga_n^{\frac{1}{n-1}} = \limsup_{n\to\infty} \ga_n^{\frac{1}{n}},
\end{equation}
we immediately obtain, by combining with \cite[Cor.~1.1.7, formula (1.13)]{StahlandTotik},

\begin{lemma}
For $\mu\in\BB$ we have
\begin{equation}\label{capacitybounds}
c_\mu \leq \liminf_{n\to\infty}\;\Cpct(K_n) \leq
\limsup_{n\to\infty}\; \Cpct(K_n) \leq \Cpct(K),
\end{equation}
where $c_\mu$ is the minimal carrier capacity.
\end{lemma}
The examples in \cite[Section 1.5]{StahlandTotik} show that all the inequalities in \eqref{capacitybounds}
can be strict. However, in this paper we only need $\liminf_{n\to\infty}\;\Cpct(K_n)>0$,
which is implied for $\mu\in\BB_+$.

Furthermore, we have
\begin{lemma}
Let $\mu\in\BB$ and choose $R>0$ so that $K\subset \D(0,R)$.
Then there exists $N$ such that for all $n\geq N$:
\begin{equation}
K_n \subset P_n^{-1}(\overline{\D(0,R)}) \subset \D(0,R).
\end{equation}
\end{lemma}
{\it Proof.}
By \cite[Theorem 1.1.4]{StahlandTotik}, we have
$$
\liminf_{n\to\infty}\frac{1}{n}\log|P_n(z)| \geq g_\Om(z)
$$
locally uniformly on $\C\sm \Co(K)$.
Taking $R$ such that $K\subset \D(0,R)$ then $\partial \D(0,R)$
is a compact set disjoint from $K$ on which $g_\Om$ is continuous, and hence
$\eps = \inf\{g_\Om(z) \,|\, |z| = R \} > 0$.
By the above inequality and compactness of $\partial\D(0,R)$,
there exists $N$ such that 
$$
\forall\; n \geq N \; \forall\; z\in \partial \D(0,R) :\quad \frac{1}{n}\log|P_n(z)| \geq \eps/2.
$$
By increasing $N$ if necessary, we can suppose $\log(R) < N\eps/2$.
Then since the zeros of $P_n$ are contained in $\Co(K) \subset \D(0,R)$ (by
Theorem \ref{thm:fejer}), the minimal modulus principle implies
$$
\forall\; n \geq N :\quad P_n(\C\sm\D(0,R)) \subset \C\sm\overline{\D(0,R)}.
$$
Thus, by \eqref{filledJuliaset},
$$
\forall\; n \geq N :\quad K_n\subset P_n^{-1}(\overline{\D(0,R)}) \subset \D(0,R)
$$
and this completes the proof.
\hfill $\square$



\begin{prop}
\label{THMlogPnvsgn}
Let $\mu\in\BB_+$. 
Then there exists $N\in\N$ and $M>0$ such that
\begin{equation}\label{logPnvsgnformula}
\forall\; n\geq N:\quad \Bigl\Vert \,g_n(z) - \frac{1}{n}\log^+|P_n(z)|\, \Bigr\Vert_\infty \leq \frac{M}{n}.
\end{equation}
\end{prop}
\begin{rem}
Proposition \ref{THMlogPnvsgn} plays a key role in the proofs of our main results. It links the Green's functions $g_n$ for $\Om_n$ to the potentials
$\frac{1}{n}\log|P_n(z)|$ or rather to the Green's functions $\frac{1}{n}\log^+|P_n(z)|$
of the set $\{z \,|\, |P_n(z)| > 1\}$.
The literature on orthogonal polynomials, and \cite{StahlandTotik}
in particular, does not seem to study the latter Green's function in connection
with orthogonal polynomials, though this restriction of $\frac{1}{n}\log|P_n(z)|$ is quite natural.
For instance, the equilibrium measure $\Laplace(\frac{1}{n}\log^+|P_n(z)|)$
on $\{z \,|\, |P_n(z)| = 1\}$ is the balayage in $\{z \,|\, |P_n(z)| < 1\}$ of the purely atomic measure
$\Laplace(\frac{1}{n}\log|P_n(z)|)$ with an atom of weight $\frac{1}{n}$ at each root of $P_n$
(counting multiplicities).
\end{rem}
{\it Proof of Proposition \ref{THMlogPnvsgn}.} By \eqref{capacityformula}--\eqref{infsup}, we have
\begin{equation}
\BB_+ = \{\mu\in\BB \,|\, \displaystyle{\liminf_{n\to\infty}}\,\Cpct(K_n) >0 \}.
\end{equation}
Hence $c:=\displaystyle{\liminf_{n\to\infty}\; \Cpct(K_n)} > 0$
and we can choose $R' >1$ such that $K\subset\D(0,R')$.
Further, let $R = 2R'$, $c' = c/2$ and choose $N$ so that 
$$
\forall\; n \geq N :\quad K_n\subset P_n^{-1}(\overline{\D(0,R')}) \subset\D(0,R')
\quad\textrm{and}\quad R'\geq \Cpct(K_n) > c'.
$$
The Green's functions $g_n$ can be written as
$$
g_n(z) = \log|z| - \log\Cpct(K_n) + \int \log|1-w/z| \,d\om_n(w),
$$
where $\om_n$ is harmonic measure from $\infty$.  Writing
$$
M' = \max\{\log R', -\log c'\} + \log 2 \; \mbox{ and } \; M = 3M',
$$
we find
$$
\forall\; n\geq N \; \forall\; z, |z| \geq R: \quad \bigl| g_n(z) - \log|z| \bigr|
< M'.
\quad
$$
For each $n$, denote by $A_n$ the set $\{z \,|\, |P_n(z)| < R\}$.
Then for each $n\geq N$ and all $z\in\C\sm A_n$, we have $|P_n(z)| \geq R$
so that $\log^+|P_n(z)| = \log|P_n(z)|$ and
$$
\left|g_n(z) - \frac{1}{n}\log^+|P_n(z)| \right| =
\left|\frac{1}{n}g_n(P_n(z)) - \frac{1}{n}\log|P_n(z)| \right| \leq \frac{M'}{n}.
$$
Moreover, for all $z\in\partial A_n$,
\begin{align*}
0 < g_n(z) &= \frac{1}{n}\log|P_n(z)|  + \biggl(g_n(z) - \frac{1}{n}\log|P_n(z)|\biggr) \\
&<  \biggl|\frac{1}{n}\log|P_n(z)| \biggr| + \left|g_n(z) - \frac{1}{n}\log|P_n(z)| \right|
< \frac{2M'}{n}.
\end{align*}
Hence, by the maximum principle for subharmonic functions,
$g_n(z) < {2M'}/{n}$ on all of $A_n$.
Since
$$
0 \leq \frac{1}{n}\log^+|P_n(z)| \leq \frac{M'}{n}
$$
on $A_n$ by construction, we have
$$
\left| g_n(z) - \frac{1}{n}\log^+|P_n(z)| \right| < \frac{3M'}{n} = \frac{M}{n}
$$
on $A_n$ and thus on all of $\C$.
\hfill $\square$

\begin{rem} (i) If $\liminf_{k\to\infty}\; \Cpct(K_{n_k})>0$ for some subsequence $\{n_k\}$,
then the proof shows that \eqref{logPnvsgnformula} holds when $n$ is replaced by
$n_k$.

(ii)
By \eqref{capacitybounds}, the hypothesis in the proposition is satisfied
if the minimal carrier capacity is strictly positive. However, according to
\cite[Example 1.5.4]{StahlandTotik},
there are measures $\mu\in\BB$ for which
$0 = c_\mu < \lim_{n\to\infty}\; \Cpct(K_{n})$.

\end{rem}

Combining Proposition \ref{THMlogPnvsgn} with  
\cite[Theorem 1.1.4]{StahlandTotik}, we can now prove the following proposition.


\begin{prop}
\label{prop:oldthmB}
For all $\mu\in\BB_+$ we have
\begin{equation}\label{newSandTupperbound}
\limsup_{n\to\infty}g_n(z) \leq g_\mu(z)
\end{equation}
locally uniformly in $\C$ and
\begin{equation}\label{newSandTlowerbound}
\liminf_{n\to\infty}g_n(z) \geq g_\Om(z)
\end{equation}
locally uniformly in $\C\setminus \Co(K)$.
In $\Co(K)\cap\Om$, the lower bound \eqref{newSandTlowerbound}
holds true only in capacity, that is,
for every compact set $V\subseteq\Om$ and every $\eps>0$, we have
\begin{equation}\label{newSandTlowerboundincapacity}
\lim_{n\to\infty}\Cpct(\{z\in V \,|\, g_n(z) < g_\Om(z) -\eps\}) = 0.
\end{equation}
\end{prop}

\begin{rem}
(i) As with \cite[(1.6)]{StahlandTotik}, the bound \eqref{newSandTupperbound}
holds for every $\mu\in\BB$.

(ii) For a sequence of real valued functions $h_n$ on an open set $U$ and $h: U\to \mathbb R$, the relation
$$
\limsup_{n\to\infty} h_n(z) \leq h(z) \ \text{locally uniformly in\ } U
$$
means that for every $z\in U$ and every sequence $\{z_n\}\subset U$ converging to $z$, we have
$\limsup_{n\to\infty} h_n(z_n) \leq h(z)$.
Similar statements hold for $\liminf$ and $\lim$.

\end{rem}
{\it Proof of Proposition \ref{prop:oldthmB}.}
If $c_\mu = 0$, then $g_\mu\equiv\infty$ and \eqref{newSandTupperbound} trivially holds.
The relations \eqref{newSandTupperbound} and \eqref{newSandTlowerbound} are straightforward translations of
the relations (1.6) and (1.7) from \cite[Theorem 1.1.4]{StahlandTotik}
by using Proposition \ref{THMlogPnvsgn} and noting that
for any $\eps >0$,
$$
\frac{1}{n}\log|P_n(z)| \leq g_\mu(z)+\eps \; \implies \; \frac{1}{n}\log^+|P_n(z)| \leq g_\mu(z)+\eps.
$$
This implication holds by definition of $\log^+$, since $g_\mu(z) \geq 0$.

For \eqref{newSandTlowerboundincapacity}, let $\eps>0$ be given and choose
according to Proposition \ref{THMlogPnvsgn} an $N$ such that
$$
\forall n\geq N \; \forall z\in \mathbb C: \quad
\left|g_n(z) - \frac{1}{n}\log^+|P_n(z)| \right| < \eps/2.
$$
Then for $n\geq N$, we have
$$
g_n(z)+\eps \geq \frac{1}{n}\log^+|P_n(z)| +\eps/2 \geq \frac{1}{n}\log|P_n(z)| + \eps/2
$$
so that $g_n(z) < g_\Om(z) - \eps$ implies $\frac{1}{n}\log|P_n(z)|  < g_\Om(z) - \eps/2$.
Hence,
$$
|P_n(z)|^{\frac{1}{n}} \leq e^{g_\Om(z) - \eps/2} = e^{g_\Om(z)} -(1 - e^{-\eps/2})e^{g_\Om(z)}
\leq e^{g_\Om(z)} -(1 - e^{-\eps/2}),
$$
recalling that $g_\Om (z)\geq 0$. Thus, with $\eps' := (1 - e^{-\eps/2}) > 0$ and $V\subset \Om$ a compact subset, we have
$$
\{z\in V \,|\, g_n(z) < g_\Om(z) - \eps \} \subseteq
\{z \in V \,|\, |P_n(z)|^{\frac{1}{n}} < e^{g_\Om(z)} - \eps' \}
$$
and \eqref{newSandTlowerboundincapacity} applies.
\hfill $\square$

%

\section{Relating the sequences $K_n$, $J_n$ to $K$ and $J$}
\label{NearlyHaussdorfconv}
\label{sec:3}
This section contains the proof of Theorem \ref{THMlimsupliminf}.
We shall equip the space of non-empty compact subsets of $\C$ with the Hausdorff distance,
which is the natural choice in dynamical systems (see, e.g., \cite{Douady}). We begin by briefly recalling the main definitions and then characterize $\liminf$ and $\limsup$ in this setup.

Let $\KK$ denote the set of non-empty compact subsets of $\C$.
For $L, M \in \KK$, we define
the Hausdorff semi-distance from $L$ to $M$ by
\begin{equation}
\dH(L, M) := \sup\{\d(z, M) \,|\, z\in L\}
= \sup_{z\in L}\;\inf_{w\in M}\; |z-w|
\end{equation}
and the Hausdorff distance between the two sets as
\begin{equation}
\DH(L, M) := \max\{\dH(L, M), \dH(M, L)\}.
\end{equation}
The Hausdorff distance is a metric on the space $\KK$ of compact subsets.
When $\{K_n\}\subset\KK$ is a bounded sequence of compact sets
(i.e., a sequence for which there exists $R>0$
such that $K_n\subset\D(0,R)$ for all $n$), we define the symbols
\begin{align}
\liminf_{n\to\infty} K_n &:= \{ z\in\C \,|\,
\exists \, \{z_n\},\, K_n\ni z_n \underset{n\to\infty}\longrightarrow z\},\\ 
\limsup_{n\to\infty} K_n &:= \{ z\in\C \,|\, \exists \, \{n_k\}, \, n_k\nearrow \infty
\textrm{ and }
\exists \, \{z_{n_k}\}, \, K_{n_k}\ni z_{n_k} \underset{k\to\infty}\longrightarrow z\}.
\end{align}
Clearly, $\liminf_{n\to\infty} K_n\subseteq \limsup_{n\to\infty} K_n$ and by Lemma \ref{liminfandlimsupcompact},
the sets
$$
I= \displaystyle{\liminf_{n\to\infty} K_n, \quad S=\limsup_{n\to\infty} K_n}
$$
are compact. 
The set $I$ may be empty whereas $S$ is always non-empty.
Moreover, to illustrate that $(\KK, \DH)$ is a nice metric space,
let us remark that it can be shown that $I$ is either empty or it
is the largest compact set for which
\begin{equation}
\displaystyle{\lim_{n\to\infty}}\dH(I,K_n) = 0.
\end{equation}
Likewise, $S$ is the smallest compact set for which
\begin{equation}
\displaystyle{\lim_{n\to\infty}}\dH(K_n,S) = 0.
\end{equation}
Thus, $I = S$ if and only the sequence $\{K_n\}$ is convergent to the common value $I = S$.
If the sequence $\{K_n\}$ is Cauchy, then the equality $I = S$ easily follows
and this  shows that $\KK$ is a complete metric space.
Also, the above statements serve to explain the names $\liminf$ and $\limsup$.
\begin{lemma}\label{liminfandlimsupcompact}
Let $\{K_n\}$ be a bounded sequence from $\KK$. The complements of $I= \displaystyle{\liminf_{n\to\infty}}\, K_n$
and $S = \displaystyle{\limsup_{n\to\infty}}\, K_n$ are open and characterized by
\begin{equation}
\label{eq:lemma-liminf}
z_0\in\C\sm I \; \Longleftrightarrow \; \exists\;\de_0>0\;
\exists\, \{n_k\},\, n_k\nearrow \infty\textrm{ s.t. }\; \forall\; k: \,
\d(z_0, K_{n_k}) > \de_0
\end{equation}
and
\begin{equation}
\label{eq:lemma-limsup}
z_0\in \C\sm S \; \Longleftrightarrow \; \exists\;\de_0>0\; \exists\; N\textrm{ s.t. }\;
\forall n\geq N: \, \d(z_0, K_n) > \de_0.
\end{equation}
As a consequence, both $I$ and $S$ are compact.
\end{lemma}
{\it Proof.}
The implication ``$\Leftarrow$'' in \eqref{eq:lemma-liminf} is trivial. For the reverse implication,
let $z_0\in\C$ and suppose the right hand side is false.
Then
$$
\forall \de>0\;\exists\; N\textrm{ s.t. }\;\forall\;n\geq N : \, \d(z_0, K_n) \leq \de.
$$
For each $n$, let $z_n\in K_n$ be a point with $|z_n-z_0| = \d(z_0, K_n)$.
Then $K_n\ni z_n\to z_0$ 
which shows that $z_0\in I$.

The implication ``$\Leftarrow$'' in \eqref{eq:lemma-limsup} is also trivial.
For the reverse implication, take an arbitrary $z_0\in\C$ and assume the right hand side is false.
Then for any $\de>0$ there are infinitely many values of $n$ for which $\d(z_0, K_n) \leq \de$.
Thus we may construct an increasing sequence
$\{n_k\}$ of integers such that $\d(z_0, K_{n_k}) \leq 1/k$, say.
Take as above, for each $k$, a point $z_k\in K_{n_k}$ with
$|z_k - z_0| = \d(z_0, K_{n_k}) \leq 1/k$.
Hence $z_0\in S$.

Openness of $\C\sm I$ and of $\C\sm S$ follow from the relations \eqref{eq:lemma-liminf} and \eqref{eq:lemma-limsup}. Thus $I$ and $S$ are both closed, and also bounded.
\hfill $\square$

\medskip


After these preliminaries we are ready to prove Theorem \ref{THMlimsupliminf}.

\bigskip

{\it Proof of Theorem \ref{THMlimsupliminf}(i).} Since $z\in K_n$ if and only if $g_n(z) = 0$ and since $g_\Om(z)>0$ on $\Om = \C\sm K$,
the inclusion
$$
\limsup_{n\to\infty} K_n \subseteq \Co(K)
$$
follows immediately from \eqref{newSandTlowerbound} and Lemma \ref{liminfandlimsupcompact}.
Next, choose $R>0$ so large that $K_n \subset \D(0,R)$ for all $n\geq 2$. For given $\eps>0$ we obtain from \eqref{newSandTlowerboundincapacity} that
$$
\lim_{n\to\infty}\Cpct(\{z\in V_\eps\cap\overline{\D(0,R)} \,|\,
g_n(z) < g_\Om(z) - \eps/2\}) = 0,
$$
where $V_{\eps}=\{z\in \C\, | \, g_{\Om}(z)\geq \eps\}$. Since $g_n(z) = 0$ on $K_n$,
we deduce that 
$$
\lim_{n\to\infty}\Cpct(V_\eps\cap K_n) = 0
$$
and the proof 
is complete.
\hfill $\square$
\bigskip

{\it Proof of Theorem \ref{THMlimsupliminf}(ii).}
Since the right hand side of the relation is closed, it suffices to prove that
$J\sm E_\mu \subseteq \displaystyle{\liminf_{n\to\infty}}\; J_n$.
Suppose to the contrary that there exists
$z_0\in J\sm E_\mu$ which does \emph{not} belong to $\displaystyle{\liminf_{n\to\infty}}\; J_n$.
Then $g_\mu(z_0) = 0$ and according to Lemma \ref{liminfandlimsupcompact},
\begin{equation}\label{Jdedisjoint}
\exists\;\de>0 \;
\exists \, \{n_k\}, \, n_k\nearrow \infty \, \textrm{ s.t. } \; \forall\; k : \,
\D(z_0, \de)\cap J_{n_k} = \emptyset.
\end{equation}
Since $z_0\in J$, there exists $w_0\in\D(z_0, \de/4)\cap\Om$.
Choose $r \leq \de/4$ such that $\overline{\D(w_0, r)}\subset\Om$.
Let $2\eps = g_\Om(w_0)>0$ and define
$$
L := \{ w \in \overline{\D(w_0, r)}\; |\; g_\Om(w) \geq 2\eps\}.
$$
Let $L_0$ denote the connected component of $L$ containing $w_0$.
Since $g_\Om$ is subharmonic, it has no local maxima.
It follows that
$L_0\subset\Om\cap\overline{\D(z_0, \de/2)}$
is a non-trivial compact continuum and hence $\Cpct(L_0)>0$.
Thus, by \eqref{newSandTlowerboundincapacity}
there exists $N$ such that
$$
\forall\;k\geq N : \; \Cpct(\{z\in L_0| g_{n_k}(z) \leq g_\Om(z) -\eps\}) < \Cpct(L_0).
$$
Since $g_\Om(z)\geq 2\eps$ on $L_0$, it follows that
$$
\forall\;k\geq N \; \exists\; z_k\in L_0 \, \textrm{ s.t. } g_{n_k}(z_k) \geq \eps.
$$
Combining with \eqref{Jdedisjoint}, we find that $\D(z_0,\de)\subset\Om_{n_k}$ for
$k\geq N$. By applying Harnacks inequality, 
we obtain
$$
g_{n_k}(z_0) \geq g_{n_k}(z_k)\frac{1-1/2}{1+1/2} \geq \eps/3 >0.
$$
On the other hand, by \eqref{newSandTupperbound},
$$
\limsup_{k\to\infty} g_{n_k}(z_0) \leq g_\mu(z_0) = 0,
$$
which is a contradiction.
\hfill $\square$

%

\section{Results for $n$-th root regular measures}
\label{sec:4}
In this section, we specialize the general results of the previous sections to the
important class of regular measures.
According to Stahl and Totik, a measure $\mu\in \BB$ is $n$th-root regular,
in short $\mu\in\Reg$, if
\begin{equation}
 \label{regdef}
\lim_{n\to\infty} {\textstyle\frac{1}{n}}\log|P_n(z)| = g_\Om(z)
\end{equation}
locally uniformly for $z\in\C\sm\Co(K)$. In particular, we see that $\Reg\subset\BB_+$.
Note that \eqref{regdef} is equivalent to \cite[Theorem~3.2.1, formula (2.1)]{StahlandTotik}
\begin{equation}
\label{SandTregbound}
\limsup_{n\to\infty}|P_n(z)|^{1/n} \leq e^{g_\Om(z)}
\end{equation}
locally uniformly in $\C$.

A prime example of $\mu\in\Reg$ is the equilibrium measure for the boundary $J$
of a full compact non-polar subset $K$ or, equivalently, the harmonic measure
on $\C\sm K$ viewed from infinity. This follows immediately from Erd{\"o}s-Tur{\'a}n's
theorem, see \cite[Theorem~4.1.1]{StahlandTotik}.

Combining \eqref{regdef}--\eqref{SandTregbound} with Proposition \ref{THMlogPnvsgn}, we have as an immediate corollary
\begin{cor}
\label{cor eq}
The following statements are equivalent:
\begin{enumerate}[(i)]
\item
$\mu\in\Reg$,
\item
$\displaystyle{\lim_{n\to\infty} g_n(z) = g_\Om(z)}$
locally uniformly for $z\in\C\sm\Co(K)$,
\item
$\displaystyle{\limsup_{n\to\infty} g_n(z) \leq g_\Om(z)}$
locally uniformly on $\C$,
\item
$\mu\in\BB_+$ and $\displaystyle{\lim_{n\to\infty} g_n(z) = 0}$ qu.\ e.\ on $J$.
\end{enumerate}
\end{cor}

Proceeding as in the proof of Theorem \ref{THMlimsupliminf}(ii), but using Corollary \ref{cor eq}(iii)
instead of \eqref{newSandTupperbound},
we obtain a stronger result (compare also with Theorem \ref{THMlimsupliminf}(i)).
\begin{cor}\label{RegJalmostinliminfofJn}
Suppose $\mu\in\Reg$. Then
\begin{equation}
\overline{J\sm E} \subseteq \liminf_{n\to\infty} J_n,
\end{equation}
where $E$ denotes the ($F_\sigma$ and polar)
exceptional set for the Green's function $g_\Om$. In particular, if $J$ is Dirichlet regular, then
\begin{equation}
\label{eq:special-case}
J\subseteq \liminf_{n\to\infty} J_n.
\end{equation}
\end{cor}
In the convex case we note the following proposition.
\begin{prop}\label{convex}
If $\mu\in\Reg$ and $K = \Co(K)$, then
\begin{equation}
J \subseteq \liminf_{n\to\infty} K_n \subseteq
\limsup _{n\to\infty} K_n \subseteq K.
\end{equation}
\end{prop}
{\it Proof.}
For a compact convex set $K$, every boundary point is Dirichlet regular.
Moreover, $J_n\subset K_n$ so that
the first inclusion follows from \eqref{eq:special-case}.
The latter follows from Theorem \ref{THMlimsupliminf}(i).
\hfill $\square$

\begin{cor}
For any compact convex subset $K$
and any $\eps>0$, there exists a polynomial $P_n$
(of high degree $n$) with
\begin{equation}
D(\partial K, K_n) < \eps \; \textrm{ and } \; D(K_n, K) < \eps.
\end{equation}
\end{cor}

It has recently been shown that a general compact connected subset $K\subset\C$
can be approximated arbitrarily well in the Hausdorff topology by (filled) Julia sets
of polynomials, see Lindsay \cite{Lindsay} and Bishop--Pilgrim \cite{BishopandPilgrim}.
Theorem \ref{THMlimsupliminf}(i), Corollary \ref{RegJalmostinliminfofJn}, and Proposition \ref{convex}
of this paper deal with approximation of general compact sets $K\subset\C$
by the (filled) Julia sets of orthogonal polynomials for probability measures
supported on $\partial K$. These results can be viewed as a complement to the results of \cite{Lindsay,BishopandPilgrim}
in the connected case and an extension in the general case.
At the same time, our results are statements about
the dynamical behaviour of orthogonal polynomials.

%

\begin{rem}
We cannot expect that
\begin{equation}
\limsup_{n\to\infty} K_n \subseteq K
\end{equation}
for general non-convex sets $K$.
To see this, suppose $K\subset \C$ is any full compact subset of $\C$ with $K = -K$
(i.e., $z\mapsto -z$ is an involution of $K$)
and let $\om$ denote the equilibrium measure on $J = \partial K$.
Then the corresponding orthonormal polynomials $P_n$ are even for $n$ even and
odd for $n$ odd.
In particular, $0$ is a fixed point of each $P_{2n+1}$, $n\geq 0$,
and so $0\in K_{2n+1}$.
This implies that $0 \in \limsup_{n\to\infty} K_n$.
However, we may choose $K$ as above with $0\notin K$.
Note that $K$ cannot be connected in this case.
\end{rem}


\section{The orthogonal polynomials for the measure of maximal entropy of a polynomial}
\label{sec:5}

Our main results apply to measures $\mu\in\BB_+$ or $\mu\in\Reg$.
A natural way of generating non-trivial examples of such measures is to take a monic polynomial $Q$
of degree $d\geq 2$ and construct the unique balanced invariant measure $\om$ for $Q$ (see, e.g., \cite{Brolin}).
This measure is known to coincide with the (unique) measure of maximal entropy for $Q$ (see \cite{Lyubich})
and is in fact the equilibrium measure of $J_Q$, the Julia set of $Q$.
Note that, with $K_Q$ the filled Julia set of $Q$, we have 
\[
\Cpct(J_Q)=\Cpct(K_Q)=1.
\] 

The orthogonal polynomials associated with $\omega$ (as above) were studied in a series of papers of Barnsley et al. One of their basic results reads:
\begin{thm}[\cite{Barnsleyetal}]
\label{BanGerHar}
Let $Q(z) = z^d + az^{d-1} + \cdots$ be a polynomial of degree $\geq 2$ and let $\om$
denote the unique measure of maximal entropy for $Q$.
Then the monic orthogonal polynomials $\{p_n\}$ with respect to $\om$ satisfy
\begin{enumerate}
\item [(i)]
$p_1(z) = z + a/d$,
\item [(ii)]
$\forall k\in\N$ : $p_{kd}(z) = p_k(Q(z))$,
\item [(iii)]
$\forall k\in\N$ : $p_{d^k}(z) = p_1(Q^k(z)) = Q^k(z) + a/d$.
\end{enumerate}
\end{thm}

The last part of this theorem in particular shows that if $Q$ is centered (i.e., $a=0$), then the iterates of $Q$ fit neatly into the sequence
of monic orthogonal polynomials. To be specific, 
\[
 Q^k=p_{d^k} \; \mbox{ for all $k\geq 0$}.
\] 
A natural question in this context is: Are the remaining orthogonal polynomials dynamically related to $Q$?
As a corollary of Theorem~\ref{THMlimsupliminf} we obtain the following answer to this question:
\begin{cor}
In the setting of Theorem \ref{BanGerHar}, let $J_n$ and $K_n$ be the Julia set, resp. filled Julia set, of the orthonormal 
polynomial $P_n=\ga_np_n$. Then
\begin{equation}
J_Q \subseteq \liminf_{n\to\infty} J_n \subseteq \limsup_{n\to\infty} K_n \subseteq \Co(K_Q).
\end{equation}
Moreover, for any $\eps>0$ and $V_\eps := \{z\in \C \,|\, g_\Om(z) \geq \eps \}$,
\begin{equation}
\lim_{n\to\infty}\Cpct(V_\eps\cap K_n) = 0.
\end{equation}
\end{cor}

{\it Proof.}
Since $\om\in\Reg$ and $J_Q$ is Dirichlet regular, this follows from Corollary \ref{RegJalmostinliminfofJn} and Theorem \ref{THMlimsupliminf}.
\hfill $\square$





              Jacob Stordal Christiansen,
              Lund University, Centre for Mathematical Sciences, Box 118, 22100 Lund, Sweden\\
              \texttt{stordal@maths.lth.se}           

      Christian Henriksen,
              DTU Compute, Technical University of Denmark, Build. 303B, 2800 Kgs. Lyngby, Denmark\\
              \texttt{chrh@dtu.dk}           

      Henrik Laurberg Pedersen,
               Department of Mathematical Sciences, University of Copenhagen,
	       Universitetsparken 5, 2100 Copenhagen, Denmark\\
              \texttt{henrikp@math.ku.dk}           

      Carsten Lunde Petersen,
              Department of Science and Environment, Roskilde University, 4000 Roskilde, Denmark \\
              \texttt{lunde@ruc.dk}           

\end{document}